\documentclass[12pt]{article}

\usepackage[margin=1in]{geometry}
\usepackage{amsthm, amssymb, amsmath, amsfonts}
\usepackage{color}

\newtheorem{theorem}{Theorem}[section]
\newtheorem{lemma}[theorem]{Lemma}
\newtheorem{proposition}[theorem]{Proposition}

\theoremstyle{definition}
\newtheorem{definition}[theorem]{Definition}
\newtheorem{example}[theorem]{Example}
\newtheorem{remark}[theorem]{Remark}

\newcommand{\mibb}[1]{\left[#1\right]}
\newcommand{\bK}{\mathbb{K}}
\newcommand{\bC}{\mathbb{C}}

\newcommand{\cA}{\mathcal{A}}
\newcommand{\cS}{\mathcal{S}}
\newcommand{\cB}{\mathcal{B}}
\newcommand{\EMR}{\mathcal{E}\mathcal{M}\mathcal{R}_\ell}
\newcommand{\pb}{\overline{\varphi}}
\newcommand{\fr}[1]{f_{[#1]}}
\newcommand{\ue}{\mathrm{e}}
\newcommand{\ff}[2]{#1^{\underline{#2}}}

\allowdisplaybreaks

\begin{document}
\sloppy

\title{Exponential Multiple Riordan Arrays}

\author{Tian-Xiao He\\
{\small Department of Mathematics}\\
{\small Illinois Wesleyan University}\\
{\small Bloomington, IL 61702-2900, USA}\\
}
\date{}
\maketitle

\begin{abstract}
\noindent
We extend exponential double Riordan arrays to exponential multiple Riordan arrays $[g;f_1,\ldots,f_\ell]$, whose columns are generated, in exponential form, by a function $g$ of $t^\ell$ and $\ell$ multipliers $f_1,\ldots,f_\ell$ used cyclically. We prove an $\ell$-case fundamental theorem, define the exponential multiple Riordan group, and determine its identity, inverses, and its Appell, associated, and $\ell$ Bell-type subgroups; the group is the semidirect product of its Appell and associated subgroups. We compute the production matrix in terms of $2\ell$ sequences, the $c_r$- and $r_r$-sequences ($0\leq r\leq \ell-1$), give its bivariate generating function in terms of generalized hyperbolic functions, and derive recurrence relations for the entries. Finally, we give two sequence characterizations: one by $\ell-1$ one-step $A$-sequences, an $\ell$-step $A$-sequence, and a $Z$-sequence, and one by a single $\ell$-step $A$-sequence and $\ell$ $Z$-sequences. For $\ell=2$ all results reduce to those for exponential double Riordan arrays.


\end{abstract}

\section{Introduction}

Riordan arrays are infinite lower triangular matrices defined by the generating functions of their columns; they form the Riordan group (Shapiro, Getu, Woan, and Woodson \cite{SGWW}). For a sequence $(c_n)_{n\geq 0}$ of nonzero numbers, the generalized Riordan array $(g,f)$ with respect to $(c_n)$ has entries
\begin{equation}\label{gen.def}
d_{n,k}=\mibb{\frac{t^n}{c_n}}g(t)\frac{f(t)^k}{c_k}=\frac{c_n}{c_k}[t^n]g(t)f(t)^k,
\end{equation}
and for $c_n=n!$ we obtain exponential Riordan arrays. The production matrix of a Riordan array $R$ is $P=R^{-1}UR$, where $U=(\delta_{i+1,j})_{i,j\geq 0}$ is the upper shift matrix (Deutsch, Ferrari, and Rinaldi \cite{DFR}; see also \cite{He15}).

Davenport, Shapiro, and Woodson \cite{DSW} introduced double Riordan arrays $(g;f_1,f_2)$ with an even function $g$ and two odd multipliers used alternately. Their sequence characterizations and compressions are studied in \cite{He18}, and the $A$- and $Z$-sequences in \cite{BDFJT, DFSW}. Triple Riordan arrays were studied by Barry \cite{Barry}, and the multiple Riordan group, with $\ell$ multipliers, was defined and characterized by one $A$-sequence and several $Z$-sequences in \cite{He-MR}. An ordinary multiple Riordan array $(g;f_1,\ldots,f_\ell)$ has columns
\[
(g,\ gf_1,\ gf_1f_2,\ \ldots,\ gf_1\cdots f_\ell,\ gf_1\cdots f_\ell f_1,\ \ldots),
\]
where $g\in\bK[[t^\ell]]$ and $f_i\in t\bK[[t^\ell]]$. In \cite{He-EDRA} we introduced exponential double Riordan arrays $[g;f_1,f_2]$ and studied the exponential double Riordan group, its production matrices, and its sequence characterizations. The purpose of this paper is to extend these results to any number $\ell\geq 2$ of multipliers. The structure of the resulting group from the point of view of fiber products is the subject of the companion paper \cite{He-FPM}.

Throughout, $\bK=\bC$ (or any field of characteristic zero in which the required $\ell$th roots exist; see \cite{He-FPM} for a root-free treatment). We use the notation
\[
\fr{0}=1,\qquad \fr{r}=f_1f_2\cdots f_r\quad (1\leq r\leq \ell),\qquad F=\fr{\ell}=f_1\cdots f_\ell,
\]
and $\ff{n}{\ell}=n(n-1)\cdots(n-\ell+1)$ for the falling factorial.

\section{Exponential multiple Riordan arrays and the exponential multiple Riordan group}

\begin{definition}\label{def:EMRA}
Let $\ell\geq 2$, and let
\[
g(t)=\sum_{k\geq 0}g_{\ell k}\frac{t^{\ell k}}{(\ell k)!},\qquad f_i(t)=\sum_{k\geq 0}f_{i,\ell k+1}\frac{t^{\ell k+1}}{(\ell k+1)!}\quad (1\leq i\leq \ell),
\]
with $g_0\not=0$ and $f_{i,1}\not=0$ for all $i$. The {\it exponential multiple Riordan array} $[g;f_1,\ldots,f_\ell]=[d_{n,j}]_{n,j\geq 0}$ is the lower triangular matrix whose columns have the exponential generating functions
\begin{equation}\label{EMRA.def}
\sum_{n\geq 0}d_{n,\ell k+r}\frac{t^n}{n!}=\frac{g\,\fr{r}\,F^k}{(\ell k+r)!},\qquad k\geq 0,\ 0\leq r\leq \ell-1,
\end{equation}
i.e., the columns are
\[
\left(g,\ \frac{gf_1}{1!},\ \frac{gf_1f_2}{2!},\ \ldots,\ \frac{g\fr{\ell-1}}{(\ell-1)!},\ \frac{gF}{\ell!},\ \frac{gFf_1}{(\ell+1)!},\ \ldots\right).
\]
\end{definition}

For $\ell=2$ this is the exponential double Riordan array of \cite{He-EDRA}. Since $g\in\bK[[t^\ell]]$ and $f_i\in t\bK[[t^\ell]]$, the column \eqref{EMRA.def} lies in $t^r\bK[[t^\ell]]$; hence
\begin{equation}\label{support}
d_{n,j}=0\qquad\mbox{unless } n\equiv j\pmod \ell.
\end{equation}
We say that a power series has {\it residue} $s$ if it lies in $t^s\bK[[t^\ell]]$.

The product $F=f_1\cdots f_\ell$ has residue $0$ and order $\ell$, with leading coefficient $\kappa=f_{1,1}\cdots f_{\ell,1}\not=0$. Fix an $\ell$th root of $\kappa$; then $F$ has a unique $\ell$th root
\[
\varphi=F^{1/\ell}=\kappa^{1/\ell}\,t\,\big(1+O(t^\ell)\big)\in t\bK[[t^\ell]],
\]
which has order $1$. We denote by $\pb$ its compositional inverse; $\pb\in t\bK[[t^\ell]]$ and $F(\pb(t))=t^\ell$. Also $f_i/\varphi\in\bK[[t^\ell]]$ has order $0$.

\begin{theorem}\label{thm:FT}
{\rm (Fundamental theorem of exponential multiple Riordan arrays)} Let $D=[g;f_1,\ldots,f_\ell]$ and $0\leq s\leq\ell-1$. If $A(t)=\sum_{k\geq 0}a_{\ell k+s}t^{\ell k+s}/(\ell k+s)!$ has residue $s$, and $A(t)$ is identified with the column vector of its coefficients $(a_0,a_1,\ldots)^T$, then
\begin{equation}\label{FT}
DA(t)=g\,\frac{\fr{s}}{\varphi^s}\,A(\varphi).
\end{equation}
For a general $A(t)=\sum_n a_nt^n/n!$, $DA$ is the sum of the right-hand sides of \eqref{FT} over the residue components of $A$.
\end{theorem}

\begin{proof}
By \eqref{EMRA.def}, the exponential generating function of $DA$ is
\[
\sum_{k\geq 0}a_{\ell k+s}\frac{g\fr{s}F^k}{(\ell k+s)!}=g\fr{s}\sum_{k\geq 0}a_{\ell k+s}\frac{\varphi^{\ell k}}{(\ell k+s)!}=g\frac{\fr{s}}{\varphi^s}\sum_{k\geq 0}a_{\ell k+s}\frac{\varphi^{\ell k+s}}{(\ell k+s)!}=g\frac{\fr{s}}{\varphi^s}A(\varphi).
\]
\end{proof}

\begin{remark}
The right-hand side of \eqref{FT} does not depend on the choice of $\kappa^{1/\ell}$: replacing $\varphi$ by $\omega\varphi$ with $\omega^\ell=1$ multiplies $A(\varphi)$ by $\omega^s$ and $\varphi^{-s}$ by $\omega^{-s}$. For $\ell=2$, $\fr{1}/\varphi=f_1/\sqrt{f_1f_2}=\sqrt{f_1/f_2}$, and Theorem \ref{thm:FT} is \cite[Theorem 2.1]{He-EDRA}.
\end{remark}

Applying Theorem \ref{thm:FT} to the columns of a second array gives the multiplication rule.

\begin{theorem}\label{thm:product}
Let $D=[g;f_1,\ldots,f_\ell]$ and $D'=[d;h_1,\ldots,h_\ell]$, and let $\varphi=F^{1/\ell}$ for $D$. Then
\begin{equation}\label{prod}
[g;f_1,\ldots,f_\ell]\,[d;h_1,\ldots,h_\ell]=\left[g\,d(\varphi);\ \frac{f_1}{\varphi}h_1(\varphi),\ \ldots,\ \frac{f_\ell}{\varphi}h_\ell(\varphi)\right].
\end{equation}
\end{theorem}

\begin{proof}
Write $G=gd(\varphi)$ and $H_i=(f_i/\varphi)h_i(\varphi)$. Then $G\in\bK[[t^\ell]]$, $H_i\in t\bK[[t^\ell]]$ has order $1$, and $H_1\cdots H_\ell=(F/\varphi^\ell)\,(h_1\cdots h_\ell)(\varphi)=(h_1\cdots h_\ell)(\varphi)$. The column $\ell k+r$ of $D'$ is $d\,h_{[r]}(h_1\cdots h_\ell)^k/(\ell k+r)!$, which has residue $r$, so by \eqref{FT} the corresponding column of $DD'$ is
\[
g\frac{\fr{r}}{\varphi^r}\,\frac{d(\varphi)h_{[r]}(\varphi)(h_1\cdots h_\ell)^k(\varphi)}{(\ell k+r)!}=\frac{G\,H_1\cdots H_r\,(H_1\cdots H_\ell)^k}{(\ell k+r)!},
\]
which is the column $\ell k+r$ of the right-hand side of \eqref{prod}.
\end{proof}

\begin{proposition}\label{prop:conj}
Let $\Lambda={\rm diag}(0!,1!,2!,\ldots)$ and let $(g;f_1,\ldots,f_\ell)$ be the ordinary multiple Riordan array of \cite{He-MR} with the same $g$ and $f_i$. Then
\begin{equation}\label{conj}
[g;f_1,\ldots,f_\ell]=\Lambda\,(g;f_1,\ldots,f_\ell)\,\Lambda^{-1},\qquad\mbox{i.e.,}\qquad d_{n,j}=\frac{n!}{j!}d^*_{n,j}.
\end{equation}
Hence $(g;f_1,\ldots,f_\ell)\mapsto[g;f_1,\ldots,f_\ell]$ is an isomorphism from the multiple Riordan group of \cite{He-MR} onto the set of exponential multiple Riordan arrays with the multiplication \eqref{prod}.
\end{proposition}

\begin{proof}
By \eqref{EMRA.def}, $d_{n,j}=\frac{n!}{j!}[t^n](\mbox{column $j$ of $(g;f_1,\ldots,f_\ell)$})$. Conjugation by $\Lambda$ preserves products.
\end{proof}

\begin{theorem}\label{thm:group}
The set of all exponential multiple Riordan arrays with the multiplication \eqref{prod} is a group, called the exponential multiple Riordan group and denoted by $\EMR$. Its identity is $[1;t,\ldots,t]$, and
\begin{equation}\label{inverse}
[g;f_1,\ldots,f_\ell]^{-1}=\left[\frac{1}{g(\pb)};\ \frac{t\pb}{f_1(\pb)},\ \ldots,\ \frac{t\pb}{f_\ell(\pb)}\right].
\end{equation}
The sets
\begin{align*}
&\cS=\{[1;f_1,\ldots,f_\ell]\},\qquad \cA=\{[g;t,\ldots,t]\},\\
&\cB_j=\{[g;f_1,\ldots,f_\ell]:\ f_j=tg\}\quad (1\leq j\leq\ell)
\end{align*}
are subgroups, called the associated, Appell, and type-$j$ Bell subgroups. The Appell subgroup is normal, and $\EMR=\cA\rtimes\cS$.
\end{theorem}

\begin{proof}
By \eqref{prod}, $\EMR$ is closed under multiplication, and \eqref{prod} is matrix multiplication, hence associative; $[1;t,\ldots,t]=I$. Let $[d;h_1,\ldots,h_\ell]$ be the right-hand side of \eqref{inverse}. Then $d(\varphi)=1/g$ since $\pb(\varphi)=t$, and
\[
\frac{f_i}{\varphi}h_i(\varphi)=\frac{f_i}{\varphi}\cdot\frac{\varphi\,\pb(\varphi)}{f_i(\pb(\varphi))}=\frac{f_i\,t}{f_i}=t,
\]
so the product is $I$. Moreover $1/g(\pb)\in\bK[[t^\ell]]$ and $t\pb/f_i(\pb)\in t\bK[[t^\ell]]$ has order $1$, so \eqref{inverse} lies in $\EMR$; since a one-sided inverse of a lower triangular matrix with nonzero diagonal is two-sided, \eqref{inverse} holds.

Closure of $\cS$ and $\cA$ is clear from \eqref{prod}. For $\cB_j$: if $f_j=tg$ and $h_j=td$, then the $j$th multiplier of the product is $(f_j/\varphi)\,\varphi\, d(\varphi)=t\,g\,d(\varphi)$, which is $t$ times the new first component. For the inverse of an element of $\cB_j$, the $j$th multiplier in \eqref{inverse} is $t\pb/(\pb g(\pb))=t/g(\pb)$, again $t$ times the first component.

Finally, the multipliers of a product depend only on the multipliers of the factors, so $\pi([g;f_1,\ldots,f_\ell])=[1;f_1,\ldots,f_\ell]$ is a homomorphism $\EMR\to\cS$ with kernel $\cA$ which is the identity on $\cS$. Since $[g;f_1,\ldots,f_\ell]=[g;t,\ldots,t][1;f_1,\ldots,f_\ell]$ and $\cA\cap\cS=\{I\}$, we get $\EMR=\cA\rtimes\cS$.
\end{proof}

We shall use the action of $D^{-1}$ on series of a given residue.

\begin{lemma}\label{lem:invaction}
Let $D=[g;f_1,\ldots,f_\ell]$ and let $B$ have residue $s$, $0\leq s\leq \ell-1$. Then
\begin{equation}\label{invaction}
D^{-1}B=\frac{t^s\,B(\pb)}{g(\pb)\,\fr{s}(\pb)}.
\end{equation}
\end{lemma}

\begin{proof}
The components of $D^{-1}$ in \eqref{inverse} are $\tilde g=1/g(\pb)$ and $\tilde f_i=t\pb/f_i(\pb)$. Since $F(\pb)=t^\ell$, $\tilde f_1\cdots\tilde f_\ell=t^\ell\pb^\ell/F(\pb)=\pb^\ell$, so the $\ell$th root for $D^{-1}$ is $\tilde\varphi=\pb$. By Theorem \ref{thm:FT},
\[
D^{-1}B=\tilde g\,\frac{\tilde f_1\cdots\tilde f_s}{\pb^s}\,B(\pb)=\frac{1}{g(\pb)}\prod_{i=1}^s\frac{t}{f_i(\pb)}\,B(\pb).
\]
\end{proof}

\section{Production matrices}

For an exponential generating function, deleting the first row of a column corresponds to differentiation. Thus $\overline D=UD$ has columns $(g\fr{r}F^k)'/(\ell k+r)!$, and the production matrix $P=(P_0,P_1,\ldots)=D^{-1}\overline D$ can be computed from Lemma \ref{lem:invaction}.

\begin{theorem}\label{thm:production}
Let $D=[g;f_1,\ldots,f_\ell]$ and define
\begin{align*}
&c_0(t)=\frac{t^{\ell-1}g'(\pb)}{g(\pb)\fr{\ell-1}(\pb)},\qquad r_0(t)=\frac{F'(\pb)}{\fr{\ell-1}(\pb)},\\
&c_r(t)=\frac{(g\fr{r})'(\pb)}{(g\fr{r-1})(\pb)},\qquad r_r(t)=f_r(\pb)\,F'(\pb),\qquad 1\leq r\leq \ell-1.
\end{align*}
Then the columns of the production matrix $P$ of $D$ have the exponential generating functions
\begin{align}
P_{\ell k}&=\frac{1}{(\ell k)!}\left(c_0(t)t^{\ell k}+k\,r_0(t)t^{\ell k-1}\right),\label{Pcol0}\\
P_{\ell k+r}&=\frac{1}{(\ell k+r)!}\left(c_r(t)t^{\ell k+r-1}+k\,r_r(t)t^{\ell k+r-1-\ell}\right),\qquad 1\leq r\leq\ell-1,\label{Pcolr}
\end{align}
and, with coefficients of negative index equal to zero,
\begin{align}
p_{n,\ell k}&=\frac{n!}{(\ell k)!}\left(c_{0,n-\ell k}+k\,r_{0,n-\ell k+1}\right),\label{p0}\\
p_{n,\ell k+r}&=\frac{n!}{(\ell k+r)!}\left(c_{r,n-\ell k-r+1}+k\,r_{r,n-\ell k-r+1+\ell}\right),\qquad 1\leq r\leq\ell-1.\label{pr}
\end{align}
\end{theorem}

\begin{proof}
The column $\ell k$ of $\overline D$ is $(gF^k)'/(\ell k)!=(g'F^k+kgF^{k-1}F')/(\ell k)!$, which has residue $\ell-1$. By \eqref{invaction} with $s=\ell-1$ and $F(\pb)=t^\ell$,
\[
P_{\ell k}=\frac{t^{\ell-1}}{(\ell k)!}\cdot\frac{g'(\pb)t^{\ell k}+k\,g(\pb)F'(\pb)t^{\ell k-\ell}}{g(\pb)\fr{\ell-1}(\pb)}=\frac{1}{(\ell k)!}\left(c_0t^{\ell k}+k\,r_0t^{\ell k-1}\right).
\]
For $1\leq r\leq\ell-1$, the column $\ell k+r$ of $\overline D$ is $\big((g\fr{r})'F^k+k\,g\fr{r}F^{k-1}F'\big)/(\ell k+r)!$, which has residue $r-1$. By \eqref{invaction} with $s=r-1$,
\[
P_{\ell k+r}=\frac{t^{r-1}}{(\ell k+r)!}\cdot\frac{(g\fr{r})'(\pb)t^{\ell k}+k\,(g\fr{r})(\pb)F'(\pb)t^{\ell k-\ell}}{(g\fr{r-1})(\pb)},
\]
which is \eqref{Pcolr}. Taking $n![t^n]$ gives \eqref{p0} and \eqref{pr}.
\end{proof}

\begin{remark}\label{rem:residues}
Since $\pb\in t\bK[[t^\ell]]$, composition with $\pb$ preserves residues. Hence $c_0\in t^{\ell-1}\bK[[t^\ell]]$, while $r_0$, $c_r$, and $r_r$ ($r\geq 1$) lie in $\bK[[t^\ell]]$; moreover $r_r$ has order at least $\ell$. Consequently,
\[
p_{n,j}=0\qquad\mbox{unless } n\equiv j-1\pmod\ell,
\]
in agreement with $P=D^{-1}UD$, where $D^{\pm1}$ satisfy \eqref{support} and $U$ shifts residues by one. Comparing leading terms gives $r_0(0)=\ell f'_\ell(0)$ and $c_r(0)=r\,f'_r(0)$ for $1\leq r\leq\ell-1$. For $\ell=2$, $(c_0,r_0,c_1,r_1)$ are the sequences $(c_1,r_1,c_2,r_2)$ of \cite[Theorem 3.1]{He-EDRA}.
\end{remark}

\subsection{The bivariate generating function}

Let $\omega=\ue^{2\pi i/\ell}$ and define the generalized hyperbolic functions
\[
E_{\ell,r}(x)=\sum_{k\geq 0}\frac{x^{\ell k+r}}{(\ell k+r)!}=\frac{1}{\ell}\sum_{j=0}^{\ell-1}\omega^{-jr}\ue^{\omega^jx},\qquad 0\leq r\leq\ell-1,
\]
with indices read modulo $\ell$. For $\ell=2$, $E_{2,0}=\cosh$ and $E_{2,1}=\sinh$. Since $E'_{\ell,r}=E_{\ell,r-1}$ and $xE'_{\ell,r}(x)=\sum_k(\ell k+r)x^{\ell k+r}/(\ell k+r)!$,
\begin{equation}\label{kE}
\sum_{k\geq 0}\frac{k\,x^{\ell k+r}}{(\ell k+r)!}=\frac{1}{\ell}\left(xE_{\ell,r-1}(x)-rE_{\ell,r}(x)\right).
\end{equation}

\begin{proposition}\label{prop:bgf}
The bivariate generating function $\phi_P(t,z)=\sum_{n,j}p_{n,j}\frac{t^n}{n!}z^j=\sum_jP_j(t)z^j$ of the production matrix is
\begin{align*}
\phi_P(t,z)=&\,c_0(t)E_{\ell,0}(zt)+\frac{z\,r_0(t)}{\ell}E_{\ell,\ell-1}(zt)\\
&+\sum_{r=1}^{\ell-1}\left(\frac{c_r(t)}{t}E_{\ell,r}(zt)+\frac{r_r(t)}{\ell\,t^{\ell+1}}\Big(zt\,E_{\ell,r-1}(zt)-r\,E_{\ell,r}(zt)\Big)\right).
\end{align*}
\end{proposition}

\begin{proof}
Multiply \eqref{Pcol0} and \eqref{Pcolr} by $z^{\ell k}$ and $z^{\ell k+r}$, sum over $k$, and use \eqref{kE} with $x=zt$. For example, $\sum_kk\,r_0t^{\ell k-1}z^{\ell k}/(\ell k)!=r_0t^{-1}\cdot\frac{1}{\ell}\,zt\,E_{\ell,\ell-1}(zt)$.
\end{proof}

For $\ell=2$ this is the corrected formula $\phi_P=\left(c_1+\frac{zr_2}{2t^2}\right)\cosh(zt)+\left(\frac{zr_1}{2}+\frac{c_2}{t}-\frac{r_2}{2t^3}\right)\sinh(zt)$ of \cite{He-EDRA}.

\subsection{Recurrence relations}

\begin{theorem}\label{thm:rec}
Let $D=[g;f_1,\ldots,f_\ell]=[d_{n,j}]$. Then, for $n\geq 0$ and $k\geq 0$ (coefficients of negative index being zero),
\begin{align*}
d_{n+1,\ell k}&=\sum_{i\geq \ell k-1}\frac{i!}{(\ell k)!}\left(c_{0,i-\ell k}+k\,r_{0,i-\ell k+1}\right)d_{n,i},\\
d_{n+1,\ell k+r}&=\sum_{i\geq \ell k+r-1}\frac{i!}{(\ell k+r)!}\left(c_{r,i-\ell k-r+1}+k\,r_{r,i-\ell k-r+1+\ell}\right)d_{n,i},\qquad 1\leq r\leq\ell-1.
\end{align*}
\end{theorem}

\begin{proof}
This is $\overline D=DP$, i.e., $d_{n+1,j}=\sum_id_{n,i}p_{i,j}$, together with \eqref{p0}, \eqref{pr}, and $p_{i,j}=0$ for $i<j-1$.
\end{proof}

\subsection{Examples}

\begin{example}\label{ex:Appell}
For an element $[g;t,\ldots,t]$ of the Appell subgroup, $\varphi=\pb=t$ and $\fr{r}=t^r$, so
\[
c_0(t)=\frac{g'(t)}{g(t)},\qquad r_0(t)=\ell,\qquad c_r(t)=r+\frac{tg'(t)}{g(t)},\qquad r_r(t)=\ell t^\ell\quad(1\leq r\leq\ell-1).
\]
\end{example}

\begin{example}\label{ex:PS}
Consider the Pascal--Stanley-type array
\[
D=\left[\frac{1}{1-t^\ell};\ t,\ \ldots,\ t,\ \frac{t}{1-t^\ell}\right],
\]
with $f_1=\cdots=f_{\ell-1}=t$ and $f_\ell=t/(1-t^\ell)$; for $\ell=2$ it is the array of \cite[Example 4.4]{He-EDRA}. Here $F=t^\ell/(1-t^\ell)$, and the column $\ell k+r$ is $t^{\ell k+r}/\big((1-t^\ell)^{k+1}(\ell k+r)!\big)$, so
\[
d_{n,\ell k+r}=\left\{\begin{array}{ll}\dfrac{n!}{(\ell k+r)!}\dbinom{(n-r)/\ell}{k},&\mbox{if } n\equiv r\pmod\ell,\\[2ex] 0,&\mbox{otherwise}.\end{array}\right.
\]
We have $\varphi=t(1-t^\ell)^{-1/\ell}$, $\pb=t(1+t^\ell)^{-1/\ell}$, and $1-\pb^{\,\ell}=1/(1+t^\ell)$. Theorem \ref{thm:production} gives
\[
c_0=\ell t^{\ell-1}(1+t^\ell),\quad r_0=\ell(1+t^\ell)^2,\quad c_r=r+\ell t^\ell,\quad r_r=\ell t^\ell(1+t^\ell)\quad(1\leq r\leq\ell-1).
\]
For $\ell=3$, the array and its production matrix begin
\[
D=\mibb{\begin{array}{cccccccccc}
1\\ 0&1\\ 0&0&1\\ 6&0&0&1\\ 0&24&0&0&1\\ 0&0&60&0&0&1\\ 720&0&0&240&0&0&1\\ 0&5040&0&0&420&0&0&1\\ 0&0&20160&0&0&672&0&0&1\\ 362880&0&0&181440&0&0&1512&0&0&1\\ \vdots&&&&&&&&&\ddots
\end{array}},
\]
\[
P=\mibb{\begin{array}{cccccccccc}
0&1\\ 0&0&1\\ 6&0&0&1\\ 0&18&0&0&1\\ 0&0&36&0&0&1\\ 360&0&0&180&0&0&1\\ 0&0&0&0&180&0&0&1\\ 0&0&0&0&0&252&0&0&1\\ 0&0&0&40320&0&0&840&0&0&1\\ 0&0&0&0&0&0&0&648&0&0\\ \vdots&&&&&&&&&\ddots
\end{array}}.
\]
For instance, $p_{2,0}=2!\,c_{0,2}=2\cdot 3=6$, $p_{3,1}=3!\,c_{1,3}=6\cdot 3=18$, and $p_{5,3}=\frac{5!}{3!}(c_{0,2}+r_{0,3})=20\cdot(3+6)=180$.
\end{example}

\section{Sequence characterizations}

By \eqref{conj}, sequence characterizations of ordinary multiple Riordan arrays translate into exponential ones. We give two such characterizations. Recall that $\varphi=F^{1/\ell}$ and $\pb$ is its compositional inverse; all generating functions below are series in $t^\ell$, and we write their coefficients as $a_j=[t^{\ell j}]A(t)$, and so on.

\subsection{One-step $A$-sequences, an $\ell$-step $A$-sequence, and a $Z$-sequence}

Let $(d^*_{n,j})=(g;f_1,\ldots,f_\ell)$ be the ordinary multiple Riordan array. Define
\begin{align}
&A_r(t)=\sum_{j\geq 0}a_{r,j}t^{\ell j}=\frac{f_r(\pb)}{\pb}\quad(1\leq r\leq\ell),\qquad A(t)=\sum_{j\geq0}a_jt^{\ell j}=\frac{t^\ell}{\pb^{\,\ell}},\label{Ar}\\
&Z(t)=\sum_{j\geq 0}z_jt^{\ell j}=\frac{1}{\pb^{\,\ell}}\left(1-\frac{g(0)}{g(\pb)}\right).\label{Z}
\end{align}
Equivalently,
\begin{equation}\label{Ar-eq}
f_r=tA_r(\varphi),\qquad F=t^\ell A(\varphi),\qquad g=\frac{g(0)}{1-t^\ell Z(\varphi)},
\end{equation}
and in particular $A_1(t)A_2(t)\cdots A_\ell(t)=A(t)$.

\begin{lemma}\label{lem:ord-I}
The ordinary multiple Riordan array $(d^*_{n,j})=(g;f_1,\ldots,f_\ell)$ satisfies
\begin{align*}
&d^*_{n,\ell k+r}=\sum_{j\geq0}a_{r,j}\,d^*_{n-1,\ell(k+j)+r-1},\qquad 1\leq r\leq\ell,\ n\geq 1,\\
&d^*_{n,\ell k}=\sum_{j\geq0}a_j\,d^*_{n-\ell,\ell(k+j-1)},\qquad k\geq 1,\ n\geq \ell,\\
&d^*_{n,0}=\sum_{j\geq0}z_j\,d^*_{n-\ell,\ell j},\qquad n\geq\ell,
\end{align*}
where in the first formula with $r=\ell$ the column $\ell k+\ell$ is meant, i.e., $d^*_{n,\ell(k+1)}=\sum_ja_{\ell,j}d^*_{n-1,\ell(k+j)+\ell-1}$.
\end{lemma}

\begin{proof}
The column $\ell k+r$ ($1\leq r\leq\ell$) equals $f_r\cdot g\fr{r-1}F^k=t\sum_ja_{r,j}\,g\fr{r-1}F^{k+j}$ by \eqref{Ar-eq}, since $A_r(\varphi)=\sum_ja_{r,j}F^j$; the right-hand side is $t$ times a combination of the columns $\ell(k+j)+r-1$. Similarly, the column $\ell k$, $k\geq1$, equals $F\cdot gF^{k-1}=t^\ell\sum_ja_jgF^{k-1+j}$, and $g-g(0)=t^\ell\,g\,Z(\varphi)=t^\ell\sum_jz_jgF^j$. Multiplication by $t$ or $t^\ell$ shifts rows by $1$ or $\ell$.
\end{proof}

\begin{theorem}\label{thm:charI}
A lower triangular matrix $[d_{n,j}]_{n,j\geq 0}$ with $d_{0,0}\not=0$ and $d_{n,0}=0$ for $1\leq n\leq\ell-1$ is an exponential multiple Riordan array if and only if there exist sequences $(a_{r,j})_{j\geq 0}$, $1\leq r\leq\ell-1$, $(a_j)_{j\geq0}$, and $(z_j)_{j\geq0}$ with $a_{r,0}\not=0$ and $a_0\not=0$ such that, for $k\geq0$,
\begin{align}
&d_{n,\ell k+r}=\frac{n}{(\ell k+r)!}\sum_{j\geq0}(\ell(k+j)+r-1)!\,a_{r,j}\,d_{n-1,\ell(k+j)+r-1},\quad 1\leq r\leq\ell-1,\ n\geq1,\label{I1}\\
&d_{n,\ell k}=\frac{\ff{n}{\ell}}{(\ell k)!}\sum_{j\geq0}(\ell(k+j-1))!\,a_j\,d_{n-\ell,\ell(k+j-1)},\quad k\geq 1,\ n\geq\ell,\label{I2}\\
&d_{n,0}=\ff{n}{\ell}\sum_{j\geq0}(\ell j)!\,z_j\,d_{n-\ell,\ell j},\quad n\geq\ell.\label{I3}
\end{align}
In this case the generating functions of the sequences are given by \eqref{Ar} and \eqref{Z}.
\end{theorem}

\begin{proof}
Necessity follows from Lemma \ref{lem:ord-I} by substituting $d^*_{n,j}=\frac{j!}{n!}d_{n,j}$ from \eqref{conj} and multiplying by $n!/j!$; for instance, $\frac{(\ell k+r)!}{n!}d_{n,\ell k+r}=\sum_ja_{r,j}\frac{(\ell(k+j)+r-1)!}{(n-1)!}d_{n-1,\ell(k+j)+r-1}$ gives \eqref{I1}.

For sufficiency, let $\varphi$ be the unique series in $t\bK[[t^\ell]]$ of order $1$ with $\varphi^\ell=t^\ell A(\varphi)$ (it exists since $a_0\not=0$), and set
\[
f_r=tA_r(\varphi)\ (1\leq r\leq\ell-1),\qquad f_\ell=\frac{\varphi^\ell}{f_1\cdots f_{\ell-1}},\qquad g=\frac{d_{0,0}}{1-t^\ell Z(\varphi)}.
\]
Then $g\in\bK[[t^\ell]]$, each $f_i\in t\bK[[t^\ell]]$ has order $1$ (as $a_{r,0}\not=0$), and $F=\varphi^\ell$. The array $[g;f_1,\ldots,f_\ell]$ satisfies \eqref{I1}--\eqref{I3} with the given sequences, and it has the same row $0$ and the same column-$0$ entries in rows $1,\ldots,\ell-1$ (namely $0$) as the given matrix. Every other entry of row $n\geq1$ is determined by \eqref{I1}--\eqref{I3} from earlier rows (columns $\ell k$, $k\geq1$, with $n<\ell$ vanish by lower triangularity), so the two matrices coincide.
\end{proof}

\begin{remark}
For $\ell=2$, Theorem \ref{thm:charI} is \cite[Theorem 4.2]{He-EDRA}, with $A_1(t)=f_1(\overline{f_{12}})/\overline{f_{12}}$ and $A_2(t)=A(t)=t^2/\overline{f_{12}}^{\,2}$. The one-step $A_\ell$-sequence of Lemma \ref{lem:ord-I} may be used instead of the $\ell$-step $A$-sequence for the columns $\ell k$, $k\geq1$:
\[
d_{n,\ell(k+1)}=\frac{n}{(\ell(k+1))!}\sum_{j\geq 0}(\ell(k+j)+\ell-1)!\,a_{\ell,j}\,d_{n-1,\ell(k+j)+\ell-1},
\]
with $A_\ell=A/(A_1\cdots A_{\ell-1})$.
\end{remark}

\subsection{An $\ell$-step $A$-sequence and $\ell$ $Z$-sequences}

Following \cite{BDFJT, DFSW} for $\ell=2$ and \cite{He-MR} for general $\ell$, one may instead use a single $\ell$-step $A$-sequence for all columns $j\geq\ell$ and one $Z$-sequence for each of the first $\ell$ columns. Let
\begin{equation}\label{kappa}
\kappa_r=d_{r,r}=g(0)f'_1(0)\cdots f'_r(0),\qquad 0\leq r\leq\ell-1,
\end{equation}
and define, with $A(t)$ as in \eqref{Ar},
\begin{equation}\label{Zr}
Z_r(t)=\sum_{j\geq0}z_{r,j}t^{\ell j}=\frac{1}{\pb^{\,\ell}}\left(1-\frac{\kappa_r\,\pb^{\,r}}{(g\fr{r})(\pb)}\right),\qquad 0\leq r\leq\ell-1.
\end{equation}
Thus $Z_0=Z$, and \eqref{Zr} is equivalent to $g\fr{r}-\kappa_rt^r=t^\ell g\fr{r}Z_r(\varphi)$.

\begin{theorem}\label{thm:charII}
A lower triangular matrix $[d_{n,j}]_{n,j\geq0}$ whose first $\ell$ rows vanish off the diagonal and whose entries $d_{r,r}$, $0\leq r\leq\ell-1$, are nonzero is an exponential multiple Riordan array if and only if there exist sequences $(a_j)$ with $a_0\not=0$ and $(z_{r,j})_{j\geq0}$, $0\leq r\leq\ell-1$, such that for all $n\geq\ell$
\begin{align}
&d_{n,j}=\frac{\ff{n}{\ell}}{j!}\sum_{i\geq0}(j+\ell(i-1))!\,a_i\,d_{n-\ell,j+\ell(i-1)},\qquad j\geq\ell,\label{II1}\\
&d_{n,r}=\frac{\ff{n}{\ell}}{r!}\sum_{i\geq0}(\ell i+r)!\,z_{r,i}\,d_{n-\ell,\ell i+r},\qquad 0\leq r\leq\ell-1.\label{II2}
\end{align}
In this case the generating functions are $A(t)$ in \eqref{Ar} and $Z_r(t)$ in \eqref{Zr}.
\end{theorem}

\begin{proof}
For the ordinary array, the column $j=\ell k+r\geq\ell$ is $g\fr{r}F^k=t^\ell\sum_ia_i\,g\fr{r}F^{k-1+i}$, which gives $d^*_{n,j}=\sum_ia_id^*_{n-\ell,j+\ell(i-1)}$; and the column $r<\ell$ satisfies $g\fr{r}-\kappa_rt^r=t^\ell\sum_iz_{r,i}\,g\fr{r}F^i$, which gives $d^*_{n,r}=\sum_iz_{r,i}d^*_{n-\ell,\ell i+r}$ for $n\geq\ell$. Substituting \eqref{conj} gives \eqref{II1} and \eqref{II2}. Conversely, let $\varphi$ be as in the proof of Theorem \ref{thm:charI} and set
\[
H_r=\frac{d_{r,r}\,t^r}{1-t^\ell Z_r(\varphi)}\quad(0\leq r\leq\ell-1),\qquad g=H_0,\]\[ f_r=\frac{H_r}{H_{r-1}}\ (1\leq r\leq\ell-1),\qquad f_\ell=\frac{\varphi^\ell}{f_1\cdots f_{\ell-1}}.
\]
Then $g\fr{r}=H_r$ and $F=\varphi^\ell$, and $[g;f_1,\ldots,f_\ell]$ satisfies \eqref{II1}--\eqref{II2} and agrees with the given matrix in the first $\ell$ rows; hence the two matrices coincide.
\end{proof}

\begin{remark}
For $\ell=2$, Theorem \ref{thm:charII} is \cite[Theorem 4.3]{He-EDRA}; in particular \eqref{Zr} with $r=1$ is the corrected $Z_2$-sequence $Z_2(t)=\frac{1}{\overline{f_{12}}^2}\big(1-\frac{g(0)f'_1(0)\overline{f_{12}}}{g(\overline{f_{12}})f_1(\overline{f_{12}})}\big)$. As noted there, $z_{r,0}$ may vanish.
\end{remark}

\begin{example}
For the Pascal--Stanley-type array of Example \ref{ex:PS} we have $\pb^{\,\ell}=t^\ell/(1+t^\ell)$, and \eqref{Ar}, \eqref{Z}, and \eqref{Zr} give
\[
A_1=\cdots=A_{\ell-1}=1,\qquad A_\ell=A=1+t^\ell,\qquad Z=Z_0=Z_1=\cdots=Z_{\ell-1}=1.
\]
Hence, by Theorem \ref{thm:charI},
\[
d_{n,\ell k+r}=\frac{n}{\ell k+r}\,d_{n-1,\ell k+r-1}\ (1\leq r\leq\ell-1),\qquad d_{n,0}=\ff{n}{\ell}\,d_{n-\ell,0},
\]
\[
d_{n,\ell k}=\ff{n}{\ell}\left(\frac{(\ell k-\ell)!}{(\ell k)!}\,d_{n-\ell,\ell k-\ell}+d_{n-\ell,\ell k}\right)\quad(k\geq1).
\]
For $\ell=3$: $d_{7,4}=\frac{7}{4}\,d_{6,3}=\frac{7}{4}\cdot240=420$ and $d_{9,3}=9\cdot8\cdot7\left(\frac{1}{6}\,d_{6,0}+d_{6,3}\right)=504\,(120+240)=181440$.
\end{example}

\section{Concluding remarks}

For $\ell=2$, all results of this paper reduce to those of \cite{He-EDRA}. The square root $\sqrt{f_1f_2}$ used there is replaced here by the $\ell$th root $\varphi=(f_1\cdots f_\ell)^{1/\ell}$, which requires $\ell$th roots of $f'_1(0)\cdots f'_\ell(0)$. In the companion paper \cite{He-FPM} we show that $\EMR$ is an $\ell$-fold fiber product of generalized Riordan groups with respect to the sequences $((\ell n+r)!)_{n\geq0}$, $0\leq r\leq\ell-1$. This gives root-free product and inverse formulas, describes the group structure, and shows that the $\ell$th power of the production matrix $P$ of Theorem \ref{thm:production} is the direct sum of the production matrices of the $\ell$ compressions of $[g;f_1,\ldots,f_\ell]$.

\end{document}